\documentclass[10pt,a4paper,oneside,leqno]{article}
\usepackage{amsfonts,amssymb}
\usepackage{amscd}
\newtheorem{defn}{Definition}

\vspace{0.5cm}
 4
\font\ebf=cmbx8
\font\erm=cmr8

\usepackage{graphicx, graphics}

\newcommand{\fnomial}[2]{ {{#1} \choose {#2}}_F }

\begin{document}
\begin{center}
	\noindent { \textsc{ How  the  work of Gian Carlo Rota had influenced my  group  research and life}}  \\ 
	\vspace{0.3cm}
	\vspace{0.3cm}
	\noindent Andrzej Krzysztof Kwa\'sniewski \\
	\vspace{0.2cm}
	\noindent {\erm Member of the Institute of Combinatorics and its Applications  }\\
{\erm High School of Mathematics and Applied Informatics} \\
	{\erm  Kamienna 17, PL-15-021 Bia\l ystok, Poland }\noindent\\
	\noindent {\erm e-mail: kwandr@gmail.com}\\
	\vspace{0.2cm}
\end{center}

\noindent {\ebf Abstract:}
\noindent {\small One outlines here in a brief overview how  the  work of Gian Carlo Rota had influenced my  research and life, starting from the end of the last century up to present time state of The Internet Gian Carlo Rota Polish Seminar.}

\vspace{0.1cm}

\vspace{0.1cm}

\noindent AMS Classification Numbers: 05A40, 81S99, 06A06 ,05B20, 05C7  

\vspace{0.1cm}

\noindent  affiliated to The Internet Gian-Carlo Polish Seminar:

\noindent \emph{http://ii.uwb.edu.pl/akk/sem/sem\_rota.htm}


\vspace{0.1cm}

\section{How did I had come over}

\vspace{0.1cm}

\noindent How and when  I did came across the  work of Gian Carlo Rota - this I do not remember. May be it was only in 1997 because
of streams of thousands of references on the so called $q$-deformations (extensions) that I was slightly involved in. Then I started to agnize, to be more fully aware  of the Gian Carlo Rota's and his friends' and disciples'  outstanding importance with  his and theirs mathematical culture main stream inherited ideas, language and goals especially there, where both analysis and combinatorics meet to enjoy the join into the alloy ore - the crystalline formation of Mathemagics.
\vspace{0.1cm}
\noindent May it be then in December or so in 1998 at Bia{\l}ystok - when  I was much impressed by a series of  Professor Oleg Viktorovich Viskov from Steklov Institute  lectures  on umbral calculus and all that.  Since that time in almost all my "`umbra"' articles I frequently refer to Professor Viskov contributions [1-4] and others - for more see  [5-7]. These [5-7] references  are examples of my first contributions (including "`upside down notation"') to the extended umbral calculus. What is this "`upside down notation"' from [5-7] ? It is just this : $k_F \equiv F_k$ , where  $F$ is a natural numbers valued sequence. This notation  inspired by Gauss  and in the spirit of Knuth via the reasoning just repeated   with   "$k_F$"  numbers   replacing  $k$ - natural numbers leads one to transparent clean results in a lot of cases as for example in the recent acyclic digraph's articles [8-10]. For this notation see also Appendix in [11].

\vspace{0.2cm}

\section{Graves-Heisenberg-Weyl algebra}

\noindent The  ingenious ideas of differential and dual graded posets that we owe to  Stanley  and  Fomin  (see [10])  
bring together combinatorics, representation theory, topology, geometry  and many more specific branches of mathematics and mathematical physics thanks to intrinsic ingredient of these mathematical descriptions which is the Graves - Heisenberg - Weyl  (GHW)  algebra usually attributed to Heisenberg by physicists and  to Herman Weyl by mathematicians  and sometimes  to both of them.

\vspace{0.1cm}

\noindent As noticed by Oleg Viktorovich Viskov in  [4]  the formula  

$$ [f(a),b] = cf'(a)$$
                                                        
\noindent where
$$[a,b]=c,  [a,c]=[b,c] = 0 $$
                                                      
\noindent pertains to Charles Graves from Dublin [12]. Then it was re-discovered by Paul Adrien Maurice Dirac and others in the next century.\\
Let us then note that the picture that emerges in [5-7]  discloses the fact that any umbral representation of finite (extended) operator calculus or equivalently - any umbral representation of GHW algebra makes up an example of the algebraization of the analysis with generalized differential operators of Markowsky  acting on the algebra of polynomials or other algebras as for example formal series algebras.

\vspace{0.5cm}

\section{Cobweb posets and DAGs named  KoDAGs}

\noindent KoDAGs are Hasse diagrams -hence directed acyclic graphs  of cobweb partially ordered sets which are  secluded in a natural way from multi-ary relations
chains' digraphs. The family of these so called cobweb posets  has been invented by the author at the dawn of this century  (for earlier references see [13,14]- for the recent ones see [8-11]) . These structures are such a generalization of the Fibonacci tree growth that allows \textbf{joint combinatorial interpretation} [13,14]  \textbf{for all of them} under the combinatorial admissibility condition.

\vspace{0.2cm}

\noindent Let $\left\{ F_n \right\}_{n\geq 0}$ be a natural numbers valued sequence with $F_0 = 1$ (or  $F_0! \equiv 0!$ being exceptional as in case of Fibonacci numbers). Any such sequence uniquely designates both $F$-nomial coefficients of an $F$-extended umbral calculus as well as $F$-cobweb poset defined in [13]. If these $F$-nomial coefficients are natural numbers or zero
then we call the sequence $F$ - the $F$-\textbf{cobweb admissible sequence}.

\vspace{0.2cm}

\begin{defn}
\noindent  Let  $n\in N \cup \left\{0\right\}\cup \left\{\infty\right\}$. Let   $r,s \in N \cup \left\{0\right\}$.  Let  $\Pi_n$ be the graded partial ordered set (poset) i.e. $\Pi_n = (\Phi_n,\leq)= ( \bigcup_{k=0}^n \Phi_k ,\leq)$ and $\left\langle \Phi_k \right\rangle_{k=0}^n$ constitutes ordered partition of $\Pi_n$. A graded poset   $\Pi_n$  with finite set of minimal 
elements is called \textbf{cobweb poset} \textsl{iff}  
$$\forall x,y \in \Phi \  i.e. \  x \in \Phi_r \ and \  y \in \Phi_s \   r \neq s\ \Rightarrow \   x\leq y   \ or \ y\leq x  , $$ 
 $\Pi_\infty \equiv \Pi. $
\end{defn}

\vspace{0.2cm}

\noindent See Fig.1.

\vspace{0.2cm}

\begin{defn}
Let any $F$-cobweb admissible sequence be given then $F$-nomial coefficients are defined as follows
$$
	\fnomial{n}{k} = \frac{n_F!}{k_F!(n-k)_F!} 
	= \frac{n_F\cdot(n-1)_F\cdot ...\cdot(n-k+1)_F}{1_F\cdot 2_F\cdot ... \cdot k_F}
	= \frac{n^{\underline{k}}_F}{k_F!}
$$
\noindent while $n,k\in \mathbb{N}$ and $0_F! = n^{\underline{0}}_F = 1$  with $n^{\underline{k}}_F \equiv \frac{n_F!}{k_F!}$ staying for falling factorial.
\end{defn}

\vspace{0.2cm}

\begin{defn}

$C_{max}(\Pi_n) \equiv  \left\{c=<x_k,x_{k+1},...,x_n>, \: x_s \in \Phi_s, \:s=k,...,n \right\} $ i.e. $C_{max}(\Pi_n)$ is the set of all maximal chains of $\Pi_n$
\end{defn}

\vspace{0.2cm}

\begin{defn}
Let  $$C_{max}\langle\Phi_k \to \Phi_n \rangle \equiv \left\{c=<x_k,x_{k+1},...,x_n>, \: x_s \in \Phi_s, \:s=k,...,n \right\}.$$
Then the $C\langle\Phi_k \to \Phi_n \rangle $ set of  Hasse sub-diagram corresponding maximal chains defines biunivoquely 
the layer $\langle\Phi_k \to \Phi_n \rangle = \bigcup_{s=k}^n\Phi_s$  as the set of maximal chains' nodes and vice versa -
for  these \textbf{graded} DAGs (KoDAGs included).
\end{defn}

\vspace{0.2cm}

\noindent The equivalent to that of [13,14] formulation of combinatorial interpretation of cobweb posets via their cover relation digraphs (Hasse diagrams) is the following.

\vspace{0.2cm}

\noindent \textbf{Theorem} \\
\noindent(Kwa\'sniewski) \textit{For $F$-cobweb admissible sequences $F$-nomial coefficient $\fnomial{n}{k}$ is the cardinality of the family of \emph{equipotent} to  $C_{max}(P_m)$ mutually disjoint maximal chains sets, all together \textbf{partitioning } the set of maximal chains  $C_{max}\langle\Phi_{k+1} \to \Phi_n \rangle$  of the layer   $\langle\Phi_{k+1} \to \Phi_n \rangle$, where $m=n-k$.}

\vspace{0.2cm} \noindent For environment needed and then  simple combinatorial proof see [14,13]  easily accessible via Arxiv.\\

\vspace{0.2cm}

\noindent 
One uses for that to proof the graded structure of Hasse diaggram and the notion  of the layer.

\vspace{0.2cm}

\vspace{0.2cm} 
\noindent \textbf{Comment 1}. For the  above Kwa\'sniewski combinatorial  interpretation of  $F$-nomials' array the diagram being directed  or not does not matter  of course, as this combinatorial interpretation is  equally valid for partitions  of the family of  $SimplePath_{max}(\Phi_k - \Phi_n)$ in  comparability graph of the Hasse  digraph with self-explanatory notation used on the way. And to this end recall: a poset is graded if and only if every connected component of its comparability \textbf{graph } is graded. We are concerned here with connected graded graphs and digraphs.

\begin{figure}[ht]
\begin{center}
	\includegraphics[width=100mm]{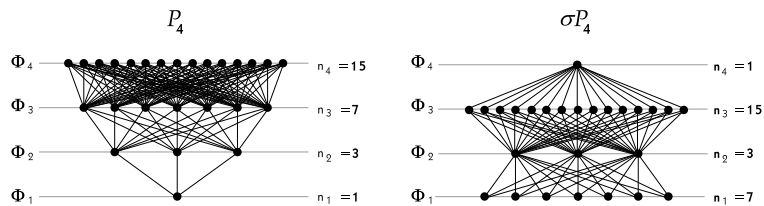}
	\caption {Display of the layer  $\langle\Phi_1 \to \Phi_4 \rangle$ = the subposet $P_4$ of the  $F$ = Gaussian integers sequence $(q=2)$ $F$-cobweb poset and $\sigma P_4$ subposet of the $\sigma$ permuted Gaussian $(q=2)$ $F$-cobweb poset .}
\end{center}
\end{figure}

\vspace{0.1cm}

\noindent \textbf{If }one  imposes further requirements with respect $F$- sequences denominating both $F$-extended Umbral (Finite Operator) Calculus and cover relation diagrams (Hasse) of the corresponding cobweb poset \textbf{then} further specific problems, their solutions and specific digraph-combinatorial interpretations are arrived at.  

\vspace{0.2cm}
\noindent For  fresh results of the Student participant of  \textbf{The Internet Gian Carlo Rota Polish Seminar} see [17].
For his  recent discoveries see  [16,17].

\vspace{0.2cm}

\section{How  all that had influenced my and my research group life ?}

\vspace{0.2cm}

\noindent The Gian Carlo Rota Polish Seminar has been transformed  in 2008 and 
is active now as The Internet Gian Carlo Rota Polish Seminar:\\
\noindent \emph{http://ii.uwb.edu.pl/akk/sem/sem\_rota.htm}. We are continuing the research.

\end{document}